\documentstyle[12pt]{article}

\newcommand{\be}{\begin{equation}}
\newcommand{\ee}{\end{equation}}
\newcommand{\bea}{\begin{eqnarray}}
\newcommand{\eea}{\end{eqnarray}}
\newcommand{\ra}{\rightarrow}

\def\Journal#1#2#3#4{{#1} {\bf #2}, #3 (#4)}

\def\PLA{{\em Phys. Lett.} A}

\def\CMP{\em Commun. Math. Phys.}
\def\LMP{\em Lett. Math. Phys.}
\def\SJNP{\em Sov. J. Nucl. Phys.}
\def\JMP{\em J. Math. Phys.}
\def\JPSJ{\em J. Phys. Soc. Jpn.}
\def\JPA{{\em J. Phys.} A}
\begin{document}

\thispagestyle{empty}

\hfill {\small DSF-T-43/98}

\bigskip

\bigskip

\begin{center}
{\Large \bf
 Tensor operators and Wigner-Eckart theorem  \\
 ~~\\
 for ${\cal U}_{q \rightarrow 0}(sl(2))$}\\

~~\\
~~\\

V. Marotta \hspace{2.5cm} A. Sciarrino\\[3mm]
Universit\`{a} di Napoli ``Federico II'' \\
Dipartimento di Scienze Fisiche \\
and \\
I.N.F.N. - Sezione di Napoli \\
I-80125 Napoli - Italy
\end{center}

\bigskip

\bigskip

\bigskip

\centerline{\bf Abstract}

\begin{quote}
Crystal tensor operators, which tranform under  ${\cal U}_{q \rightarrow
0}(sl(2))$, in analogous way as the vectors of the crystal basis, are
introduced. The Wigner-Eckart theorem for the crystal tensor is defined:
the selection rules depend on the initial state and on the component of the
tensor operator; the transition amplitudes to the states of the same final
irreducible representation are all equal.
\end{quote}

\vfill
\noindent {\small \bf
\begin{tabular}{l}
~$\,$Postal adress: Mostra d'Oltremare, Pad. 20 - I-80125 Napoli (Italy)\\
\begin{tabular}{ll}
E-mail: VMAROTTA (SCIARRINO)@NA.INFN.IT \end{tabular}
 \end{tabular} }

\pagebreak
\setcounter{page}{1}

\setcounter{equation}{0}

Deformation of enveloping Lie algebra   ${\cal U}_q({\cal G})$
introduced by Drinfled-Jimbo \cite{D}, \cite{J} is by now a subject of standard
text book. For the arguments discussed in this paper see \cite{BL}
where an accurate list of references can be found.
In the limit $q \rightarrow 0$ it has been shown by M. Kashiwara \cite{Kashi}
that  ${\cal U}_q({\cal G})$ admits a canonical
 peculiar basis, called {\it crystal basis}. Since that article crystal
 bases have been object of
very intensive mathematical studies and have also been extended
to the case of deformation of affine Kac-Moody algebras. However
a point is still , to our knowledge, missing: it is possible
to introduce the concept of q-tensor and $q$-Wigner-Eckart theorem in the
limit $q \rightarrow 0$ ?
Besides the mathematical interest, the question may be interesting in
application in physics. It is clear that in this limit we are no
more dealing with the deformation of an univeral enveloping Lie
algebra, but it is interesting to study what are the relics of
the symmetry structure described originally by the algebra ${\cal G}$
and then by the deformation of its enveloping algebra ${\cal U}_q({\cal G})$.
 It is, indeed,  well know that Wigner-Eckart theorem is one of the milestones
in the application of algebraic methods in physics. Let us remark that one of the
motivations to  study  the limit  $q \rightarrow 0$  by
Date, Jimbo and Miwa \cite{DJM}, which firstly discovered the peculiar
behavior of n-dimensional ${\cal U}_{q \rightarrow 0} (gl(n,{\cal C}))$-modules,
whose axiomatic  settlement has been given in
\cite{Kashi},  was the study of solvable lattice models where the
parameter $q$ plays the role of the temperature. Moreover
 in  \cite{FSS} the quantum enveloping algebra
 ${\cal U}_{q}(sl(2) \oplus sl(2))$ in the limit $ q \ra  0 $
 has been proposed as symmetry algebra for the genetic code assigning
 the (4) nucleotides (elementary constituents  of the genetic code) to the
  fundamental representation and the (64) codons (triplets of nucleotides )
  to
 the three-fold tensor product of the fundamental representation, using
 crystal basis.

 In the following we will consider only the crystal  basis for
 ${\cal U}_q({sl(2)})$.

To set the notation, let us recall the definition of ${\cal U}_q(sl(2))$
\be
~[J_{+},\, J_{-}]  =  [2J_{3}]_{q}
\ee
\be
~[J_3,\, J_{\pm}]  =  {\pm} \, J_{\pm}
\ee
where

\be
~[x]_q = \frac{q^x - q^{-x}} {q - q^{-1}}
\label{eq: 1}
\ee
In the following we shall omit the lower label $q$

For later use let us remind that:
\be
[n]_q! = [1]_q \: [2]_q \ldots [n]_q \label{eq: fa}
\ee
The algebra $U_{q}(sl(2))$ is endowed with an Hopf structure. In
particular we recall the coproduct is defined by
\bea
\Delta (J_{3}) & =  & J_3 \otimes {\bf 1} + {\bf 1} \otimes J_3
 \nonumber \\
\Delta (J_{\pm}) & = & J_{\pm} \otimes q^{J_3} +
q^{-J_3} \otimes J_{\pm} \label{eq: 2}
\eea

 The Casimir operator can be written
 \be
 C = J_{+}  J_{-} \, + \, [J_{3}][J_{3} - 1] =
    J_{-}  J_{+} \, + \, [J_{3}][J_{3} + 1]
 \ee

 For $q$ generic, i.e. not a root of unity, the irreducible representations
 (IR) are lalelled by
 an integer or half-integer number $j$ and the action of the generators
 on the vector basis $|jm>$, ($-j \leq m \leq j$) , of the IR is
\be
 J_{3} \,|jm> = m \, |jm>     \label{eq: 3}
 \ee
 \be
 J_{\pm} \, |jm> = \sqrt{[j \mp m] [j {\pm} m + 1]} \, |j,m {\pm} 1>
 = F^{\pm}(j,m) \,  |j,m {\pm} 1>  \label{eq: 4}
 \ee
 From eqs.(\ref{eq: 3})-(\ref{eq: 4}) it follows
 \be
 C \, |jm> = [j] [j + 1] \, |jm>
 \ee
 Let us study the behavior of a $q$-number $[x]$ for $ q \ra 0$.
 In the following the symbol $ \sim $ in the equations
 has to be read equal in the limit $q \ra 0$ modulo the addition of a function
 regular in $q = 0$.
 From the definition eq.(\ref{eq: 1})  we have
\be
  [x]_{q \ra 0}\, \sim \,  q^{-x + 1} \;\;\;\;\;\;\;\;\; \mbox{ $x \neq  0$}
 \label{eq: 5}
 \ee
 So it follows that
 \be
  F^{\pm}(j,m)_{q \ra 0} \, \sim \,  q^{-j + 1/2} \label{eq: 6}
  \ee
   \be
  [j] [j + 1]_{q \ra 0} \, \sim \,  q^{-2j + 1} \label{eq: 7}
 \ee
  \be
  [x]!_{q \ra 0}  \, \sim \,  q^{-1/2 \, x \, (x - 1)}
 \ee

 From eqs.(\ref{eq: 4}) and (\ref{eq: 6}) it follows that the action of the
 generator $J_{\pm}$ is not defined in the limit $q \ra 0$.
 Let us define the element $\Gamma_{0}$ belonging to the center of
 ${\cal U}_{q}(sl(2))$
 \be
  \Gamma_{0} = C^{-1/2} \label{eq: D}
 \ee
  \be
 \Gamma_{0} \, |jm> =  ([j] [j + 1])^{-1/2} \, |jm>_{q \ra 0}  \sim \,
 q^{j - 1/2} \, |jm>
 \ee
 Let us define
 \be
 \widetilde{J}_{\pm}  = \Gamma_{0} \, J_{\pm}
 \ee
 These operators are well behaved for ${q \ra  0}$. Their action
  in the limit $q \ra  0$ will define the crystal basis:
\bea
\widetilde{J}_{+} \, |jm> & = & |j,m+1> \quad \mbox{for} \,\, -j \leq  m < j \\
\widetilde{J}_{-} \,  |jm> & = & |j,m-1>  \quad \mbox{for} \,\, -j < m \leq  j
\eea
\be
\widetilde{J}_{+} \, |jj>  = \widetilde{J}_{-} \,  |j,-j>  =  0
\ee

The tensor product of two representations in the crystal basis
is given by \cite{Kashi}.

{\bf Theorem} - If ${\cal B}_{1}$ and ${\cal B}_{2}$
are the crystal bases of the $M_{1}$ and $M_{2}$ ${\cal U}_{q
\rightarrow 0} (sl(2))$-modules, for $u \in {\cal B}_{1}$ and $v \in
{\cal B}_{2}$, we have:
\bea
&& \tilde J_{-}(u \otimes v) = \left\{
\begin{array}{ll}
\tilde J_{-}u \otimes v & \exists \, n \ge 1 \mbox{ such that }
\tilde J_{-}^nu \ne 0 \mbox{ and } \tilde J_{+}^nv = 0 \\
u \otimes \tilde J_{-}v & \mbox{otherwise} \\
\end{array} \right. \\
&& \tilde J_{+}(u \otimes v) = \left\{
\begin{array}{ll}
u \otimes \tilde J_{+}v & \exists \, n \ge 1 \mbox{ such that }
\tilde J_{+}^nv \ne 0 \mbox{ and } \tilde J_{-}^nu = 0 \\
\tilde J_{+}u \otimes v & \mbox{otherwise} \\
\end{array} \right.
\eea

So the tensor product of two crystal basis is a crystal basis and
the  states of the basis of the tensor space are pure states. In other words in
the limit  $q \ra 0$ all the $q$-Clebsch-Gordan  ($q$-CG) coefficients vanish
except one which is equal to ${\pm} 1$..
Let us recall the definiton of $q$-tensor for  $U_{q}(sl(2))$
\cite{BT}, \cite{Nomu}, \cite{RS} and \cite{BL}.
An irreducible  $q$-tensor of rank $j$ is a family of $2j + 1$
operators $T_{m}^{j}$ ($ -j \leq m \leq j$) which tranform under
the action of the generators of  $U_{q}(sl(2))$  as
\be
q^{J_{3}}(T_{m}^{j}) \equiv q^{J_{3}} \, T_{m}^{j} \, q^{-J_{3}}
= q^{m} \, T_{m}^{j} \label{eq: C}
\ee
or
\be
~[J_{3},\, T_{m}^{j} ]  = m \,  T_{m}^{j}
\ee
\be
J_{\pm}( T_{m}^{j}) \equiv J_{\pm}\, T_{m}^{j} \, q^{J_{3}} \, - \, q^{-J_{3} {\pm} 1}
\,  T_{m}^{j} \, J_{\pm} = F^{\pm}(j,m) \, T_{m {\pm} 1}^{j}
\label{eq: CO}
\ee
In deriving the above equations use has been made of the non trivial
coproduct eq.(\ref{eq: 2}).
The q-Wigner-Eckart ($q$-WE) theorem now reads \cite{STK}
\be
<JM | T_{m}^{j} | j_{1} m_{1}> = (-1)^{2j} \, \frac{<J || T^{j} || j_{1}>}
{\sqrt{[2J + 1]}} \, <j_{1} m_{1} j m | JM> \label{eq: WE}
\ee
where  $<J|| T^{j}|| j_{1}>$ is the reduced matrix element of the
q-tensor $T^{j}$ and

 $<j_{1} m_{1}j m | JM >$ is the $q$-CG
coefficients. In the following we will use the explicit expression
of the $q$-CG of \cite{STK}.
It is useful rewrite the $q$-WE theorem eq.(\ref{eq: WE})
in the following form
\be
 T_{m}^{j} \, |j_{1} m_{1}> = (-1)^{2j} \, \sum_{J = |j - j_{1}|}^{j + j_{1}}
  \, \frac{<J || T^{j} || j_{1}>}{\sqrt{[2J + 1]}}
  \, <j_{1} m_{1} jm | JM> \, |JM> \label{eq: WES}
\ee
Our strategy to define the ($q \ra 0$)-tensor and than the
($q \ra 0$)-WE is the following:

\begin{enumerate}

\item Let us write eq.(\ref{eq: CO}) in the form
\be
 J_{\pm}\, T_{m}^{j} \, q^{J_{3}} =
q^{-J_{3} {\pm} 1} \,  T_{m}^{j} \, J_{\pm} \, + \,  F^{\pm}(j,m) \, T_{m {\pm} 1}^{j}
  \label{eq: COS}
\ee
than we multiply both sides of eq.(\ref{eq: COS}) from left and right
by an element (not unique) $\Gamma$ of the center of the algebra and define
\be
\widehat {T}^{j}_{m} = \Gamma \, T^{j}_{m} \, \Gamma
\ee
Let us remark that $\widehat{T}^{j}_{m}$ is still a
$q$-tensor operator of the same rank as $T^{j}_{m}$. Indeed it
tranforms  under the action of  $J_{{\pm}, 3}$ according to eqs.(\ref{eq: C}),
(\ref{eq: CO}) or (\ref{eq: COS}) which
have been derived by application of the coproduct eq.(\ref{eq: 2}).
 We make the {\em conjecture} that an element $\Gamma$ exists such that
 $\widehat{T}^{j}_{m}$ has a smooth and
defined behaviour in the limit $q \ra 0$.
 We will discuss below some explicit examples in which $ T^{j}_{m}$
 is not defined in the limit $q \ra 0$  and its reduced
matrix element diverge, while on the contrary  it is possible to define
 $\widehat{T}^{j}_{m}$ with a well defined limit.

\item We apply the $\widetilde{J}_{\pm}, J_{3}$ generators to eq.(\ref{eq: WES})
written for $\widehat{T}^{j}_{m}$ for  $j = 1/2$ and then we study the limit
 $q \ra 0$, assuming that
$<J || \widehat{T}^{j} || j_{1}>$ has a well-defined behaviour in the limit.

\item From the study of 2) we deduce the action of the generators
$\widetilde{J}_{\pm}, J_{3}$ in the limit $q \ra 0$ on $\widehat{T}^{1/2}_{m}$

\item From the tensor product  we can infer the action for the
generic tensor..

\end{enumerate}

\noindent To perform our second step we need to compute the $q \ra 0$ limit
of $<j_{1} m_{1} \frac{1}{2} m | JM> $. The result are reported in Table 1,
where we have used the expressions of

 $<j_{1} m_{1} \frac{1}{2} m | JM> $
given in App.B of \cite{STK} and eq.(\ref{eq: 5}).

\begin{table}[htbp]
\begin{center}
\begin{tabular}{|c|c|c|}
\hline
J &  m = 1/2 & m = -1/2 \\ \hline
$j_{1} + 1/2$ & $q^{j_{1} - m_{1}}$ & 1 \\
\hline
$j_{1} - 1/2$ & -1 & $q^{j_{1} - m_{1}}$ \, q \\
 \hline
\end{tabular}
\end{center}
\caption{Behaviour of the q-CG  $<j_{1} m_{1} \frac{1}{2} m | JM> $
 for $q \ra 0$.\label{table1}}
\end{table}

Using the results of Table 1, denoting by $\tau_{m}^{1/2}$ the
q-tensor operator $\widehat{T}_{m}^{1/2}$  in the limit $q \ra 0$ and
\be
 \left( \frac{<J || \widehat{T}^{j} || j_{1}>}
{\sqrt{[2J + 1]}} \right)_{q \ra 0} \!\!\! \sim \, <J  || \tau^{j} ||
j_{1}>
\ee
we get
\bea
 \tau_{1/2}^{1/2} \, |j_{1} m_{1}> & =  & (-1)\,
  \delta_{j_{1},m_{1}} \,<j_{1} + 1/2 \, || \tau^{1/2} || j_{1}> \,
 |j_{1} + 1/2, m_{1} + 1/2>  \nonumber \\
  & +   & <j_{1} - 1/2 \, || \tau^{1/2} || j_{1}>
 \,  |j_{1} - 1/2, m_{1} + 1/2>  \label{eq: S1}
\eea
\be
 \tau_{-1/2}^{1/2} \, |j_{1} m_{1}> = (-1) \,
  <j_{1} + 1/2 \, || \tau^{1/2} || j_{1}> \, |j_{1} + 1/2, m_{1} - 1/2>
 \label{eq: S2}
\ee
Inspection of eq.(\ref{eq: S1}) and eq\'{E}(\ref{eq: S2}) shows that the
r.h.s. of the equations  has the structure of the tensor product
of $ \underline{\frac{1}{2}} \, \otimes \, \underline {j}$
in the crystal basis, see the above quoted Kashiwara's Theorem.
Note that the {\bf order of the tensor product} is important.

So we can write the action of the generators $\widetilde{J}_{\pm}, J_{3}$ on
$\tau_{m}^{1/2}$, as:
\be
J_{3}(\tau_{m}^{1/2}) \equiv
 m \, \tau_{m}^{1/2} \;\;\;\;\;\;
\widetilde{J}_{\pm} \, ( \tau_{m}^{1/2}) \equiv  \tau_{m {\pm} 1}^{1/2} \label{eq: QT}
\ee
Clearly, if $ |m| > 1/2 $ then $\tau_{m}^{1/2}$  has to be considered vanishing.
Eq.(\ref{eq: QT})  gives for the transformation of the tensor operator
the same law as for the  crystal basis.
It has been proven by Rittenberg-Scheunert \cite{RS} that for quasitriangular
Hopf algebra ($U_{q}(sl(2))$ is "almost" quasitriangular which does not
affect the following considerations) the tensor product of tensor operators is
a tensor operator. So by applying the Rittenberg-Scheunert's theorem and in the
limit $q \ra 0$ the Kashiwara's theorem  we can extend eq.(\ref{eq: QT})
to any value $j$. So we define {\em  crystal tensor} of rank $j$ a set
of operator which transform under  $\widetilde{J}_{\pm}, J_{3}$  according to
eq.(\ref{eq: QT}). As an explicit check and a further example, we compute
the ($q \ra 0$)-WE theorem for $T^{1}$.
We need to compute the $q \ra 0$ limit
of $<j_{1} m_{1} 1 m | JM> $. The result are reported in Table 2,
where we have used the expressions of $<j_{1} m_{1} 1 m | JM> $
given in App.B of \cite{STK} and eq.(\ref{eq: 5}).

\begin{table}[htbp]
\begin{center}
 \begin{tabular}{|c|c|c|c|}
\hline
J &  m = 1 & m =  0 & m = -1 \\ \hline
$j_{1} + 1$ & $q^{2(j_{1} - m_{1})}$ & $q^{j_{1} - m_{1}}$ & 1 \\
\hline
$ j_{1} $ & $-q^{-1} q^{j_{1} - m_{1}}$ &
$q^{2}q^{2(j_{1} - m_{1})} - 1 + \delta_{j_{1}m_{1}}$ & $q q^{j_{1} - m_{1}}$ \\
\hline
$j_{1} - 1$ & 1 & $-q^{j_{1} - m_{1}}$  & $q q^{2(j_{1} - m_{1})}$\\
 \hline
\end{tabular}
\end{center}
\caption{Behaviour of the q-CG  $<j_{1} m_{1} 1 m | JM> $
 for $q \ra 0$. \label{table2}}
\end{table}

Using the results of Table 2, we obtain in the limit $ q \ra 0$:
\bea
 \tau_{1}^{1} \, |j_{1} m_{1}> & =  & <J = j_{1} + 1 \, || \tau^{1} || j_{1}>
  \,  |J, m_{1} + 1> \;\;\;\;\;\; \mbox{if} \;\;\;\;\;\; m_{1} = j_{1}
  \nonumber \\
 & =  & - <J = j_{1}  \, \, || \tau^{1} || j_{1}>  \,  |J, m_{1} + 1>
 \;\;\;\;\;\; \mbox{if} \;\;\;\;\;\; m_{1} = j_{1} - 1 \nonumber \\
 & = &  <J = j_{1} - 1 \, || \tau^{1} || j_{1}> \,  |J, m_{1} + 1>
 \;\;\;\;\;\; \mbox{if} \;\;\;\;\;\; m_{1} <  j_{1} - 1 \label{eq:  V1}
\eea
\bea
 \tau_{0}^{1} \, |j_{1} m_{1}> & = & <J = j_{1} + 1  \, || \tau^{1} || j_{1}>
  \,  |J, m_{1} >
 \;\;\;\;\;\; \mbox{if} \;\;\;\;\;\; m_{1} = j_{1}  \nonumber \\
 & = & - <J = j_{1} \, \, || \tau^{1} || j_{1}>
  \,  |J, m_{1} > \;\;\;\;\;\; \mbox{if} \;\;\;\;\;\; m_{1} < j_{1}
  \label{eq:  V2}
\eea
\be
 \tau_{-1}^{1} \, |j_{1} m_{1}> \,  =   <j_{1} + 1 \, || \tau^{1} || j_{1}>
 \,  |J, m_{1} - 1>  \label{eq:  V3}
\ee
Let us now proof the following statement:
\newtheorem{guess}{Proposition}
\begin{guess}
 If the $q$-tensors $ \widehat{T}^{r_{1}}$ and
$\widehat{T}^{r_{2}}$  have a well defined behaviour for $q \ra 0$,
i.e. the crystal tensors $\tau^{r_{1}}$  and  $\tau^{r_{2}}$ are defined,
than the $q$-tensors $\widehat{T}^{R}$, obtained by the tensor product of
  $ \widehat{T}^{r_{1}}$ and $\widehat{T}^{r_{2}}$ has a well defined
  limit for $q \ra 0$.
\end{guess}

{\bf Proof}: Let us define
\be
\widehat{T}_{K}^{R} =  \sum_{k_{1}, k_{2}} \, <r_{1} k_{1} r_{2} k_{2} | R K>
 \widehat{T}_{k_{1}}^{r_{1}} \, \widehat{T}_{k_{2}}^{j_{2}}
\label{eq: TP}
\ee
 Take the matrix element of the r.h.s. and l.h.s. of eq.(\ref{eq: TP})
 between the initial state $|j_{1} m_{1}>$ and the final state $|J M>$.
 Insert the identity
 \be
 {\bf\Large 1} =\sum_{j, m} \, |jm><jm|  \label{eq: ID}
 \ee
 in the r.h.s. and apply the $q$-WE theorem eq.(\ref{eq: WE}) for $ \widehat{T}^{r_{1}}$ and
$\widehat{T}^{r_{2}}$ . We get
\bea
<JM| \,\widehat{T}_{K}^{R} \,|j_{1} m_{1}>\!\!\! &=& \!\!\!\sum_{k_{1},
k_{2}, j, m}
  <r_{1} k_{1} r_{2} k_{2} | R K> \,  <j m r_{1} k_{1} | J M > \,
 <j_{1} m_{1} r_{2} k_{2} |jm >  \nonumber \\
&{\times}& \frac{<J \, || \,\widehat{T}^{r_{1}} \, || \,j>}{\sqrt{2J + 1}}
\; \;
 \frac{<j \, || \, \widehat{T}^{j_{2}}\, ||\,j_{1}>}{\sqrt{2j + 1}}
 \label{eq: PRO}
 \eea
  If we apply the  $q$-WE to the l.h.s. of the above equation and
  make the limit $q \ra 0$, as  by assumption the r.h.s. of eq.(\ref{eq: PRO})
  has a limit, it follows that
 \be
 \left(\frac{<J || \widehat{T}^{R} || j_{1}>}
{\sqrt{[2J + 1]}} \right)_{q \ra 0} \!\!\! \sim \, <J  || \tau^{j} ||
j_{1}>
\ee
 Use of eq.(\ref{eq: ID}) requires at least a comment. The completeness of
 the basis $|jm>$ for $su(2)$ is a particular case of the completeness
 of the IRs of a compact group. For $q \neq 1$ we cannot appeal to
 this general property as we are no more dealing with a Lie group.
 However the completeness of the $q$-coherent states \cite{G} for the
 $q$-bosons \cite{B}, \cite{M} gives us an argument for the completeness of the
 states $|jm>$ , as a realization of
 the deformed enveloping algebra $U_{q}(su(2))$, for $q$ generic, and of its
 representations can be written in terms of $q$-bosons. See below for comments
 about the use of $q$-bosons in the $q \ra 0$ limit.
   Let us remark that the knowledge of the elements $\Gamma_{r_{1}}$
and  $\Gamma_{r{2}}$, which allow to define respectively the crystal
tensors $\tau^{r_{1}}$  and  $\tau^{r_{2}}$  from $q$-tensors $T^{r_{1}}$  and
 $T^{r_{2}}$,  does not determine the element $\Gamma_{R}$, which
allows to define the crystal tensor $\tau^{R}$  from the  $q$-tensor
$T^{R}$,
obtained from the tensor product of  $T^{r_{1}}$  and  $T^{r_{2}}$,
as the elements of the center   of the algebra
 do not commute with the  generic $q$-tensors $T^{j}$, for $j \neq 0$.

Now let us discuss in some explict examples our {\em conjecture} that
it is posssible to find an element $\Gamma$ in the center of the
algebra such that the operator $ \Gamma \, T^{j} \, \Gamma $
is well defined in the limit $q \ra 0$.

\begin{itemize}

\item  Let us consider the vector operator constructed
with the generators \cite{STK}
\be
T_{\pm}^{1} = {\pm} \, \frac{1}{\sqrt{[2]}} \, q^{-J_{3}} \, J_{\pm}
\ee
\bea
T_{0}^{1}  & = & \frac{1}{[2]} \, ( q^{-1}[2J_{3}] \, + \, (q -
q^{-1}) \, J_{+}J_{-}) \nonumber \\
& = & \frac{1}{[2]} \, ( q^{-1}[2J_{3}] \, + \, (q -  q^{-1}) \,
(C - [J_{3} - 1/2]^{2})) \label{eq: VO}
\eea
has no well defined meaning in the limit $q \ra 0$. If we compute
the reduced matrix element $< j_{1}  || T^{1} || j_{1}>  $ (which
is the only non vanishing) for the $q$-tensor eq.(\ref{eq: VO}) from
eq.(\ref{eq: WE}) we get
\be
< j_{1}  || T^{1} || j_{1}> = \frac{\sqrt{[2j_{1}][2j_{1} + 1][2j_{1} + 2]}}
{[2]} \label{eq: RMV}
\ee
 and
\be
 < j_{1}  || T^{1} || j_{1}>_{q \ra 0} \sim \, q^{-3j_{1} + 1}
 \ee
 Then, choosing $\Gamma = \sqrt{q^{1/2} \, \Gamma_{0}^{3}}$, we have from
  eqs.(\ref{eq: D})-(\ref{eq: 7})
 \be
 < j_{1}  || \widehat{T}^{1} || j_{1}>_{q \ra 0} \sim \, 1
 \ee
 Let us remark that the moltiplication of $\Gamma$ by a real number, the
 addition of any element of the center vanishing for $q \ra 0$ as well
 as any functional construction of $\Gamma_{0} $ behaving in the limit
 $q \ra 0$ as  $q^{-3j_{1} + 1}$ does
 not modify our conclusion. Our choice is the {\em  minimal} one.

\item The $U_{q}(sl(2))$ can be realized in terms of $q$-bosons
\cite{B}, \cite{M} which have no well defined behaviour in the limit
as it can be seen from the defining expression
\be
a_i a_j^{+} - q^{ \delta_{ij}} \, a_j^{+} a_i = \delta_{ij} q^{-N_i}
\ee
\be
[N_i, \, a_j^{+}] = \delta_{ij} a_j^{+} ~~~~~~~ [N_i, \, a_j] = -
\delta_{ij} a_j ~~~~~ [N_i,\, N_j] = 0
\ee
or from the relation
between $q$-bosons and standard bosonic operators \cite{S}.
 So we cannot extend the $q$-boson realization to the limit $q \ra 0$.
Using $q$-boson $q$-spinor operator have been constructed \cite{BT}.
\be
T^{\frac{1}{2}}_{\frac{1}{2}} = a^{+}_{1} \, q^{\frac{N_{2}}{2}}
\;\;\;\;\;\;\; T^{\frac{1}{2}}_{-\frac{1}{2}} = a^{+}_{2} \,
q^{-\frac{N_{1}}{2}}
\label{eq: SO}
\ee
 Howevere it is always possible to compute the $q$-spinor
 reduced matrix using eq.(\ref{eq: SO}), even if in the limit $q \ra 0$
 the explicit realization of the tensor operator in terms of
 $q$-bosons is meaningless.
Indeed in the case of the $q$-vector operator above defined,
obviously we obtain the same result using the definition
 eq.(\ref{eq: VO}) in terms of the abstract generators of $U_{q}(sl(2))$
 or making use of the explicit realization of the algebra in terms of the
  $q$-bosons.
 From eq.(\ref{eq: WE}) and the expression of q-CG we get
\be
< j_{1} + 1/2 \, || T^{\frac{1}{2}} || j_{1}> = -
 \sqrt{[2j_{1} + 1][2j_{1} + 2]}
\ee
 and
\be
 < j_{1} + 1/2 \,|| T^{\frac{1}{2}} || j_{1}>_{q \ra 0} \sim \,
 q^{-2(j_{1} + 1/4)}
 \ee
 Then, choosing $\Gamma = \sqrt{q} \, \Gamma_{0}$, we have from
  eqs.(\ref{eq: D})-(\ref{eq: 7})
 \be
 < j_{1} + 1/2 \, || \widehat{T}^{\frac{1}{2}} || j_{1}>_{q \ra 0} \sim \, -1
 \ee
  If we use  $q$-spinor operator hermitean conjugate to operator given
  in eq.(\ref{eq: SO}), \cite{BL}, which is
\be
T^{\dag \,, \frac{1}{2}}_{\frac{1}{2}} = -a_{2} \, q^{-\frac{(N_{1} +
1)}{2}}
\;\;\;\;\;\;\; T^{\dag \,,\frac{1}{2}}_{-\frac{1}{2}} = a_{1} \,
q^{\frac{(N_{2} + 1)}{2}}
\label{eq: SOH}
\ee
 From eq.(\ref{eq: WE}) and the expression of $q$-CG we get
\be
< j_{1} - 1/2 \, || T^{\dag \, , \frac{1}{2}} || j_{1}> =
- \sqrt{[2j_{1}][2j_{1} + 1]}
\ee
 and
\be
 < j_{1} - 1/2 \,|| T^{\dag \, , \frac{1}{2}} || j_{1}>_{q \ra 0} \sim \,
  q^{-2(j_{1} - 1/4)}
 \ee
 Then, choosing $\Gamma = \sqrt{q} \, \Gamma_{0}$, i.e. the same value as
 for $T^{\frac{1}{2}}$  we have from eqs.(\ref{eq: D})-(\ref{eq: 7})
 \be
 < j_{1} - 1/2 \, || \widehat{T}^{\dag \, , \frac{1}{2}} || j_{1}>_{q \ra 0}
 \sim \, -1
 \ee

\end{itemize}
In conclusion we have introduced ($q \ra 0$)-tensor operators, which we
call crystal tensor operators,
making the conjecture that such operators can be obtained as the
limit for $q \ra 0$  of the $q$-tensor operator multiplied to the
right and to the left by an element $\Gamma$ of the center of the algebra.
Let us emphasize that the choice of the element  $\Gamma$  is not
unique.  We have in some specific examples shown that our conjecture is realized
 and we have explicitly determined a (minimal up a factor) form of
  $\Gamma$.
The transformation law for the generic crystal tensor operators is
\be
J_{3}(\tau_{m}^{j}) \equiv
 m \, \tau_{m}^{j} \;\;\;\;\;\;
\widetilde{J}_{\pm} \, ( \tau_{m}^{j}) \equiv  \tau_{m {\pm} 1}^{j}
\ee
Clearly, if $ |m| > j $ then $\tau_{m}^{j}$  has to be considered vanishing.
The ($q \ra 0$)-Wigner-Eckart theorem can be written
\bea
 \tau_{m}^{j} \, |j_{1} m_{1}> & =  & (-1)^{2j} \, \sum_{\alpha = 0}^{2j}
  < j_{1} + j + \alpha || \tau^{j} || j_{1}>  \,
  | j_{1} + j + \alpha,  m_{1} + m> \, \nonumber \\
 &  &  (\delta_{m_{1}, j_{1} - \alpha}
  + \delta_{-m, j - \alpha} -
   \delta_{m_{1}, j_{1} - \alpha} \, \delta_{m, j - \alpha})
   \label{eq: CT}
\eea
Let us stress that while the $q$-WE theorem (for $q$ generic)
has the same form of the usual WE theorem, roughly speaking one has
to replace the numerical expression by q-numerical expression, so
its content (selection rules, relation between the transition
amplitudes) is of the same form, the ($q \ra 0$)-WE theorem has a completely
different structure. The final IR depends not only from the rank
of the tensor and initial IR, but in a crucial way from the initial
state and from the component of the tensor in consideration.
In Table 3 we report the selection rules for the case of a
$(q \ra 0)$-vector operator, for for $\j_{1} = \frac{1}{2}, 1, \frac{3}{2} $.

\begin{table}[htbp]
\begin{center}
\begin{tabular}{|c||c|c|c|}
\hline
$m_{1} / m$  & 1 &  0 & -1 \\ \hline
$\frac{1}{2}$ & $ \frac{3}{2}$ &  $ \frac{3}{2}$  & $ \frac{3}{2}$  \\
\hline
$ - \frac{1}{2}$  & $ \frac{1}{2}$  & $ \frac{1}{2}$  & $ \frac{3}{2}$   \\
\hline  \hline
1 & 2 & 2  & 2 \\ \hline
0 & 1 & 1  & 2 \\ \hline
-1 & 0 & 1  & 2 \\ \hline  \hline
$\frac{3}{2}$  & $ \frac{5}{2}$  & $ \frac{5}{2}$  & $ \frac{5}{2}$   \\
\hline
 $\frac{1}{2}$ & $ \frac{3}{2}$ &  $ \frac{3}{2}$  & $ \frac{5}{2}$  \\
\hline
$ - \frac{1}{2}$  & $ \frac{1}{2}$  & $ \frac{3}{2}$  & $ \frac{5}{2}$   \\
\hline
$ - \frac{3}{2}$  & $ \frac{1}{2}$  & $ \frac{3}{2}$  & $ \frac{5}{2}$   \\
\hline
\end{tabular}
\end{center}
\caption{Selection rules for vector operator $\tau^{1}$. In the entries
$m_{1}/m$ the value of final $J$ in function of the components of $\tau^{1}$,
and of the initial state, for $\j_{1} = \frac{1}{2}, 1, \frac{3}{2} $.
  \label{table3}}
\end{table}

In particular the highest weight state of the initial IR $j_{1}$ is
always transformed under action of $ \tau^{j}$ into a state of
the final IR $J = j_{1} + j$, while the lowest weight state is
transformed into a state of any final IR (exactly one state if
 $j_{1} \geq j$ with $J = j_{1} + m$).
 Let us remark the peculiar feature that no vector crystal operator can be build
up with the generators $\widetilde{J}_{\pm}, J_{3}$. Indeed such a
vector crystal operator should connect any initial state to a state
of the same IR and this is not the case as one can realize from the
Table 3 or from the general form of the theorem eq.(\ref{eq: CT})
 The transitions between an initial state, belonging to IR
 $j_{1}$, and any final state, belonging to
 the IR $J$, are all equal.

\end{document}